\documentclass[11pt,a4paper]{amsart}
\usepackage[T1]{fontenc}
\usepackage[utf8]{inputenc}
\usepackage{lmodern,microtype}
\usepackage{amssymb,mathtools}
\usepackage[margin=3cm]{geometry}
\usepackage[hidelinks]{hyperref}

\newtheorem{theorem}{Theorem}[section]
\newtheorem{lemma}[theorem]{Lemma}
\newtheorem{proposition}[theorem]{Proposition}
\newtheorem{mainthm}{Theorem}

\theoremstyle{definition}
\newtheorem{definition}[theorem]{Definition}
\theoremstyle{remark}
\newtheorem{remark}[theorem]{Remark}
\numberwithin{equation}{section}

\newcommand{\R}{\mathbb{R}}
\newcommand{\C}{\mathbb{C}}
\newcommand{\Z}{\mathbb{Z}}

\newcommand{\id}{\mathrm{id}}
\DeclareMathOperator{\im}{Im}
\DeclareMathOperator{\re}{Re}
\newcommand{\abs}[1]{\lvert #1\rvert}
\newcommand{\mean}[1]{[#1]}
\newcommand{\ca}{c_{\mathrm A}}
\newcommand{\cR}{c_{\mathrm R}}
\newcommand{\cs}{c_{\ast}}
\newcommand{\Cs}{C_{\ast}}
\newcommand{\eA}{\varepsilon_{\mathrm A}}
\newcommand{\pb}[2]{\{#1,#2\}}

\title[$C^1$ invariant tori for $C^{2n}$ perturbations]{Differentiable invariant tori in nearly integrable Hamiltonian systems with $C^{2n}$ perturbations\\ satisfying a Dini condition}
\author{Joachim Paulusch}
\thanks{Earlier publications of the author appeared under the name Joachim Albrecht.}

\begin{document}

\begin{abstract}
We consider a Hamiltonian system with $n\ge2$ degrees of freedom which is a small perturbation of an analytic, integrable and non-degenerate system. The frequency vector $\omega$ of the unperturbed torus is assumed to satisfy a Diophantine condition with the smallest possible exponent $n-1$. The perturbation is $2n$ times continuously differentiable, and the moduli of continuity of its pure derivatives of order $2n$ satisfy a Dini condition, which is as such weaker than any H\"older condition. We prove that the perturbed system has an invariant torus with frequency $\omega$, which is a continuously differentiable Lagrangian graph, and on which the flow is $C^1$ conjugate to the linear flow with frequency $\omega$. In an earlier version of this result the torus was only shown to be continuous. The proof follows Moser's method of analytic smoothing; the analytic KAM theorem used in each step is proved in full, based on R\"ussmann's optimal estimate for the linearized equation.
\end{abstract}

\maketitle

\section{Introduction}

\subsection{The problem}
We consider Hamiltonian systems
\begin{equation}\label{eq:hamsys}
  \dot x=H_y(x,y),\qquad \dot y=-H_x(x,y),
\end{equation}
with $n\ge2$ degrees of freedom, where $x\in\R^n$ are angular variables (that is, $H$ is $2\pi$-periodic in $x_1,\dots,x_n$) and $y\in\R^n$ are action variables. We assume that $H$ is a small perturbation of an analytic integrable Hamiltonian $\tilde H(y)$. Expanding $\tilde H$ around $y=0$ and collecting the perturbation and the terms of order three and higher in a function $h$, we may write
\begin{equation}\label{eq:H}
  H(x,y)=a+\langle\omega,y\rangle+\tfrac12\langle y,Cy\rangle+h(x,y),
\end{equation}
with $a\in\R$, the frequency vector $\omega\in\R^n$ and a real symmetric matrix $C$; here $\langle x,y\rangle=\sum_jx_jy_j$. If $h=0$, the torus $\{y=0\}$ is invariant, and the flow on it is the linear flow $x(t)=x_0+\omega t$. KAM theory asks whether this torus survives a small perturbation $h$.

By the theorem of Kolmogorov \cite{Kolmogorov54} and Arnold \cite{Arnold63} (see \cite{Moser62} for the corresponding result for twist maps), it does, provided $H$ is analytic, $h$ is small, $C$ is invertible (Kolmogorov's non-degeneracy condition), and $\omega$ satisfies a \emph{Diophantine condition}
\begin{equation}\label{eq:dioph}
  \abs{\langle k,\omega\rangle}\ge\frac{\gamma}{\abs k^{\tau}}\qquad\text{for all }k\in\Z^n\setminus\{0\},
\end{equation}
with constants $\gamma>0$ and $\tau\ge n-1$; here and below $\abs k=\max_j\abs{k_j}$. For $\tau>n-1$ almost every $\omega$ satisfies \eqref{eq:dioph} with some $\gamma>0$. For $\tau=n-1$, which is the case of this paper, these $\omega$ form a set of measure zero, which nevertheless has the cardinality of the continuum in every open set \cite{Russmann75}. For introductions to the theory we refer to \cite{Arnold89,Moser73,Russmann79}.

\subsection{How much smoothness is needed?}
Moser's first theorem on invariant curves of area preserving twist maps \cite{Moser62} required $333$ derivatives. This was lowered by R\"ussmann \cite{Russmann70,Russmann83} and Herman \cite{Herman83}; Herman also showed that less than three derivatives are not enough. For Hamiltonian systems with $n$ degrees of freedom, Moser \cite{Moser70} outlined a proof for $h\in C^\ell$, $\ell>2n+2$, by approximating $h$ by analytic functions and applying an analytic KAM theorem to one correction at a time. P\"oschel \cite{Poschel80} carried this out for $\ell>2n$; see also \cite{Zehnder75,Zehnder76,Poschel82} for other approaches. Salamon \cite{Salamon04} proved the persistence of a single torus for $\ell>2\tau+2> 2n$, and Bounemoura \cite{Bounemoura20} and Koudjinan \cite{Koudjinan20} proved the persistence of a set of tori of positive measure under similar assumptions. In the other direction, Cheng and Wang \cite{ChengWang13} showed that a Lagrangian torus with a given frequency of an integrable system can be destroyed by perturbations which are arbitrarily small in $C^{2n-\varepsilon}$. Whether $C^{2n}$ alone suffices is, in general, an open question \cite{Wang23}.

For $\tau=n-1$, the author showed in \cite{AlbrechtDiss,Albrecht07} that the H\"older condition beyond $C^{2n}$ can be replaced by a Dini condition on the derivatives of order $2n$. There, however, the invariant torus was only obtained as a continuous manifold. Tong and Li \cite{TongLi25} later proved, for $\tau>n-1$ and unperturbed systems that need be neither analytic nor integrable, results in terms of general moduli of continuity, including the differentiability of the torus. For vector fields on the torus, the analogous problem with $n$ derivatives was treated successfully in \cite{Albrecht05F,Paulusch26}, and P\"oschel \cite{Poschel21} showed that a condition of $L^2$ type suffices there.

The present article applies the proof of \cite{Albrecht05F,Paulusch26} to Hamiltonian systems.

\subsection{Result}
For a continuous function $g\colon\R^{2n}\to\R$ which is $2\pi$-periodic in all variables let $K_g(\delta)=\sup\{\abs{g(z)-g(z')}:\abs{z-z'}\le\delta\}$ denote its modulus of continuity.

\begin{mainthm}\label{thm:main}
Let $n\ge2$, let $H$ be given by \eqref{eq:H} with $a\in\R$, a real symmetric invertible matrix $C$, and $\omega$ satisfying \eqref{eq:dioph} with $\tau=n-1$. Let $h\colon\R^{2n}\to\R$ be $2\pi$-periodic in all $2n$ variables and $2n$ times continuously differentiable, and suppose that the moduli of continuity $K_j$ of $\partial^{2n}h/\partial z_j^{2n}$, $j=1,\dots,2n$, satisfy
\begin{equation}\label{eq:dini}
  \int_0^1\frac{K_j(\delta)}{\delta}\,d\delta<\infty .
\end{equation}
Then there is a constant $\eA>0$, depending only on $n$, $\gamma$, $C$ and $K_1+\dots+K_{2n}$, such that the following holds if $\sup_{\R^{2n}}\abs h\le\eA$. There are continuously differentiable maps $X,Y\colon\R^n\to\R^n$ such that $X-\id$ and $Y$ are $2\pi$-periodic, $X$ is a $C^1$ diffeomorphism of $\R^n$, and:
\begin{enumerate}
\item For every $\xi_0\in\R^n$ the function $t\mapsto\bigl(X(\xi_0+\omega t),Y(\xi_0+\omega t)\bigr)$ is a solution of \eqref{eq:hamsys}.
\item There is a function $\Phi\in C^2(\R^n)$ with $2\pi$-periodic gradient such that $\nabla\Phi=Y\circ X^{-1}$ and $H(x,\nabla\Phi(x))$ is constant. The torus $\{(x,\nabla\Phi(x)):x\in\R^n\}$ is invariant under \eqref{eq:hamsys}, every orbit on it is dense in it, and it is a $C^1$ submanifold.
\end{enumerate}
\end{mainthm}

In other words, the invariant torus is a $C^1$ Lagrangian graph: it is the graph of the gradient of a function, so the symplectic form $\sum_j dy_j\wedge dx_j$ vanishes on it.

The periodicity of $h$ in $y$ is a technical assumption, which can be removed by a cut-off \cite[Lemma~1.6]{Poschel80}. Every H\"older condition on the derivatives of order $2n$ implies \eqref{eq:dini}, but not conversely; a modulus of continuity which behaves like $(\log1/\delta)^{-a}$ with $a>1$ is an example.

Theorem~\ref{thm:main} improves the main result of \cite{AlbrechtDiss,Albrecht07}, where $X$ was only shown to be a homeomorphism and $\Phi$ to be of class $C^1$. Compared with \cite{TongLi25}, it covers the case $\tau=n-1$, but only for analytic integrable unperturbed systems.

\subsection{Idea of the proof}
We follow Moser's method \cite{Moser70}, as suggested by R\"ussmann.\footnote{Talk by H.~R\"ussmann at the IH\'ES, Bures-sur-Yvette, on 8 March 1979.} The perturbation $h$ is approximated by analytic functions $h_k$ on complex domains of width $\delta_k\to0$, and the analytic KAM theorem is applied to one correction $h_{k+1}-h_k$ at a time. The resulting transformations $Z_k$ are composed. The key observation is the same as for vector fields \cite{Paulusch26}: the correction $h_{k+1}-h_k$ has the size $K(\delta_k)\,\delta_k^{2n}$, and the analytic theorem loses exactly the factor $\delta_k^{2n}=\delta_k^{2\tau+2}$, both in the derivative of the transformation and, after multiplication with $\delta_k$, in the transformation itself. Hence the derivatives of the transformations differ from the identity by terms of order $K(\delta_k)$, and their distances from the identity on the torus are of order $K(\delta_k)\,\delta_k$. Since $\sum_kK(\delta_k)$ is comparable to the integral in \eqref{eq:dini}, the maps $Z_k$ converge in $C^1$ on the torus.

The reason why only continuity was obtained in \cite{AlbrechtDiss,Albrecht07} is that the analytic theorem used there contained an estimate for the derivative of the transformation, but none for the transformation itself. The convergence of the $Z_k$ then had to be obtained from the Arzel\`a--Ascoli theorem. The missing estimate is \eqref{eq:Wzero} below; it follows from the proof of the analytic theorem without additional effort.

\subsection{Organization}
Section~\ref{sec:tools} contains the notation and auxiliary results. Section~\ref{sec:analytic} contains the analytic KAM theorem and its complete proof. Theorem~\ref{thm:main} is proved in Section~\ref{sec:proof}.

\section{Notation and auxiliary results}\label{sec:tools}

\subsection{Norms and domains}
For $z\in\C^m$ we use the maximum norm $\abs z=\max_j\abs{z_j}$, and for a matrix $A\in\C^{m\times l}$ the corresponding operator norm $\abs A=\max_i\sum_j\abs{a_{ij}}$. Then $\abs{Az}\le\abs A\abs z$, $\abs{AB}\le\abs A\abs B$, $\abs{A^T}\le m\abs A$ for $A\in\C^{m\times l}$, and $\abs{\langle z,w\rangle}\le m\abs z\abs w$ for $z,w\in\C^m$. The identity matrix is $E$. If $A$ is invertible and $\abs{A^{-1}(B-A)}\le\frac12$, then $B$ is invertible and $\abs{B^{-1}}\le2\abs{A^{-1}}$ (Neumann series).

For $r,s>0$ we define the complex domains
\[
  S(r)=\{x\in\C^n:\abs{\im x}<r\},\qquad D(r,s)=\{(x,y)\in\C^{2n}:\abs{\im x}<r,\ \abs y<s\},
\]
and $S'(\rho)=\{z\in\C^{2n}:\abs{\im z}<\rho\}$. All of them are convex. We denote by $P(r,s)$ the space of bounded analytic functions $F\colon D(r,s)\to\C$ which are $2\pi$-periodic in $x_1,\dots,x_n$ and real for real arguments, and by $P(r)$ the corresponding space of functions on $S(r)$. Vector and matrix valued functions of this kind form the spaces $P_m(r,s)$, $P_m(r)$, $P_{n\times n}(r)$. The supremum of $\abs F$ over a set $M$ is written $\abs F_M$, and we abbreviate $\abs F_{r,s}:=\abs F_{D(r,s)}$, $\abs F_r:=\abs F_{S(r)}$. For $f\in P_m(r)$, $\mean f=(2\pi)^{-n}\int_{[0,2\pi]^n}f(x)\,dx$ denotes the mean value.

All vectors are column vectors, and we identify $z=(x,y)\in\C^{2n}$ with the column vector $\left(\begin{smallmatrix}x\\y\end{smallmatrix}\right)$. For a vector valued function $F$ we denote by $DF$ its Jacobian matrix; in particular $V_x=(\partial V_i/\partial x_j)$ for $V(x)\in\C^n$, and $W_\zeta$ for a map $W(\zeta)\in\C^{2n}$. For a scalar function $F(x,y)$, however, $F_x$ and $F_y$ denote the gradients, which are column vectors, $\nabla F=\left(\begin{smallmatrix}F_x\\F_y\end{smallmatrix}\right)$, and $F_{yy}$ denotes the Hessian with respect to $y$; thus $DF=(\nabla F)^T$. We write $\partial_\omega f=\sum_j\omega_j\,\partial f/\partial x_j$, applied componentwise to vector valued $f$.

\subsection{Two elementary estimates}
\begin{lemma}\label{lem:cauchy}
Let $G\subseteq\C^m$ be open and $F\colon G\to\C^l$ analytic and bounded by $B$.
\begin{enumerate}
\item \emph{(Cauchy's estimate)} If the ball $\{z':\abs{z'-z}<\sigma\}$ lies in $G$, then $\abs{DF(z)}\le B/\sigma$.
\item \emph{(Mean value inequality)} If $G$ is convex and $\abs{DF}\le L$ on $G$, then $\abs{F(z)-F(z')}\le L\abs{z-z'}$ for $z,z'\in G$.
\end{enumerate}
\end{lemma}

\begin{proof}
(1) For $w\in\C^m$ with $\abs w\le1$, the function $t\mapsto F(z+tw)$ is analytic for $\abs t<\sigma$ and bounded by $B$, so $\abs{DF(z)w}\le B/\sigma$ by Cauchy's inequality in one variable; take the supremum over $w$. (2) follows from $F(z)-F(z')=\int_0^1DF(z'+t(z-z'))(z-z')\,dt$.
\end{proof}

Applying (1) to $x\mapsto F(x,y)$ or to $y\mapsto F(x,y)$, we get for example: if $F\in P(r,s)$ and $0<\delta<r$, then $\abs{F_x}_{r-\delta,s}\le\abs F_{r,s}/\delta$, $\abs{F_y}_{r,s/2}\le2\abs F_{r,s}/s$, and $\abs{F_y(x,0)}\le\abs F_{r,s}/s$ for $x\in S(r)$.

We also need the following estimate for the Taylor remainder of order three. For an analytic function $f$ on $\{\eta\in\C^n:\abs\eta<\sigma\}$ put
\[
  T_3f(\eta):=f(\eta)-f(0)-\langle f_\eta(0),\eta\rangle-\tfrac12\langle\eta,f_{\eta\eta}(0)\eta\rangle .
\]

\begin{lemma}\label{lem:taylor}
If $\abs f\le B$ on $\{\abs\eta<\sigma\}$ and $0<\varepsilon<1$, then $\abs{T_3f(\eta)}\le\dfrac{B}{1-\varepsilon}\,\dfrac{\abs\eta^3}{\sigma^3}$ for $\abs\eta\le\varepsilon\sigma$.
\end{lemma}

\begin{proof}
Fix $\eta\ne0$ with $\abs\eta\le\varepsilon\sigma$. The function $g(t)=f(t\eta)=\sum_{k\ge0}g_kt^k$ is analytic for $\abs t<\sigma/\abs\eta$ and bounded by $B$, so $\abs{g_k}\le B(\abs\eta/\sigma)^k$. Since $T_3f(\eta)=\sum_{k\ge3}g_k$, we get $\abs{T_3f(\eta)}\le B(\abs\eta/\sigma)^3\sum_{j\ge0}\varepsilon^j$.
\end{proof}

\subsection{Simple canonical transformations}
Let $J=\left(\begin{smallmatrix}0&E\\-E&0\end{smallmatrix}\right)\in\R^{2n\times2n}$. Then \eqref{eq:hamsys} reads $\dot z=J\nabla H(z)$ with $z=(x,y)$. A continuously differentiable map $W\colon G\to\C^{2n}$, $G\subseteq\C^{2n}$ open, is called \emph{canonical} if $W_\zeta^TJW_\zeta=J$ on $G$. This is equivalent to $W_\zeta JW_\zeta^T=J$, and it implies that $W$ maps solutions of the Hamiltonian system with Hamiltonian $H\circ W$ to solutions of the system with Hamiltonian $H$: if $\dot\zeta=J\nabla(H\circ W)(\zeta)=JW_\zeta^T\nabla H(W(\zeta))$, then $\frac{d}{dt}W(\zeta)=W_\zeta JW_\zeta^T\nabla H(W(\zeta))=J\nabla H(W(\zeta))$. The time-$t$ maps of Hamiltonian flows are canonical \cite{Arnold89}.

\begin{definition}
Let $r,s>0$. An analytic canonical map $W=(u,w)\colon D(r,s)\to\C^{2n}$ is called a \emph{simple canonical transformation} if its first component $u$ depends only on $\xi$, and if $W-\id\in P_{2n}(r,s)$.
\end{definition}

The composition of simple canonical transformations is again one, provided it is defined. The following lemma is \cite[Satz~A.17]{AlbrechtDiss}.

\begin{lemma}\label{lem:simple}
Let $W=(u,w)\colon D(r,s)\to\C^{2n}$ be a simple canonical transformation. Then $u_\xi(\xi)$ is invertible,
\begin{equation}\label{eq:simple}
  w(\xi,\eta)=w(\xi,0)+u_\xi(\xi)^{-T}\eta\qquad\text{and}\qquad u_\xi^T\,w_\xi=w_\xi^T\,u_\xi\quad\text{on }D(r,s).
\end{equation}
\end{lemma}

\begin{proof}
Writing out $W_\zeta^TJW_\zeta=J$ with $W_\zeta=\left(\begin{smallmatrix}u_\xi&0\\w_\xi&w_\eta\end{smallmatrix}\right)$ gives $u_\xi^Tw_\xi-w_\xi^Tu_\xi=0$ and $u_\xi^Tw_\eta=E$. Hence $u_\xi$ is invertible and $w_\eta=u_\xi^{-T}$ does not depend on $\eta$, so $w$ is affine in $\eta$.
\end{proof}

By \eqref{eq:simple}, a simple canonical transformation extends to an analytic map on $S(r)\times\C^n$, given by the same formula. By the identity theorem, the extension is canonical; it is periodic in $\xi$ and real for real arguments, but in general not bounded.

\subsection{R\"ussmann's estimate}
The small divisors enter through the following estimate for the equation $\partial_\omega u=g$, due to R\"ussmann \cite[Theorem~1.1]{Russmann75}; see also \cite[Satz~9.7]{Russmann79} and \cite[Satz~5.3]{AlbrechtDiss}. The power $\delta^{-\tau}$ is optimal.

\begin{lemma}\label{lem:russmann}
Let $\omega$ satisfy \eqref{eq:dioph} with $\tau\ge n-1$, and let $r>0$. For every $g\in P(r)$ with $\mean g=0$ there is exactly one analytic, $2\pi$-periodic function $u$ on $S(r)$ with $\mean u=0$ and $\partial_\omega u=g$; it is real for real arguments, and
\[
  \abs u_{r-\delta}\le\frac{\cR}{\gamma\,\delta^{\tau}}\,\abs g_r\qquad\text{for all }\delta\in(0,r),
\]
where $\cR>0$ depends only on $n$ and $\tau$.
\end{lemma}

\subsection{Analytic approximation}
The following theorem is \cite[Satz~2.23]{AlbrechtDiss}, specialized to functions that are periodic in all variables; compare \cite{Moser70} and Jackson's theorem \cite{Achieser67}.

\begin{theorem}[Approximation]\label{thm:approx}
Let $\nu,m\ge1$ and let $f\colon\R^\nu\to\R$ be $2\pi$-periodic in each variable and $m$ times continuously differentiable. Let $K$ be an increasing function which majorizes the moduli of continuity of $\partial^mf/\partial z_j^m$, $j=1,\dots,\nu$. Let $(\varrho_k)_{k\ge0}$ be a decreasing null sequence of positive numbers with $\varrho_0\le1$. Then there are a constant $\ca>0$, depending only on $\nu$ and $m$, and functions $f_k$ which are analytic and bounded on $S'(\varrho_k)=\{z\in\C^\nu:\abs{\im z}<\varrho_k\}$, $2\pi$-periodic in each variable and real for real arguments, such that
\begin{align*}
  \abs{f-f_k}_{\R^\nu}&\le\ca\,K(\varrho_k)\,\varrho_k^{\,m}, &
  \abs{\nabla f-\nabla f_k}_{\R^\nu}&\le\ca\,K(\varrho_k)\,\varrho_k^{\,m-1},\\
  \abs{f_k-f_{k-1}}_{S'(\varrho_k)}&\le\ca\,K(\varrho_{k-1})\,\varrho_{k-1}^{\,m}, &
  \abs{f_k}_{S'(\varrho_k)}&\le\ca\,\abs f_{\R^\nu}.
\end{align*}
\end{theorem}

\section{The analytic KAM theorem}\label{sec:analytic}

\subsection{Statement}
Throughout this section, $n\ge2$, $\tau\ge n-1$, $\gamma>0$, $\omega$ satisfies \eqref{eq:dioph}, and $C\in\R^{n\times n}$ is symmetric and invertible. We say that a function $N\in P(r,s)$ is a \emph{normal form} if there is $a\in\R$ with
\[
  N(x,0)=a\quad\text{and}\quad N_y(x,0)=\omega\qquad\text{for all }x\in S(r).
\]
Then $\{y=0\}$ is an invariant torus of the Hamiltonian system with Hamiltonian $N$, on which the flow is $x(t)=x_0+\omega t$.

\begin{theorem}[Analytic KAM theorem]\label{thm:analytic}
There are constants $\cs\in(0,1]$ and $\Cs\ge1$, depending only on $n$, $\tau$, $\gamma$, $\abs C$ and $\abs{C^{-1}}$, with the following property. Let $0<r\le1$, $s:=r^{\tau+1}$, $a\in\R$, let $Q\in P_{n\times n}(r)$ be symmetric with
\begin{equation}\label{eq:Qhyp}
  \abs{Q-C}_r\le\frac{1}{4\abs{C^{-1}}},
\end{equation}
and let $R\in P(r,s)$ with $\abs R_{r,s}\le\vartheta s^2$, where $0<\vartheta\le\cs$. Consider
\[
  H(x,y)=a+\langle\omega,y\rangle+\tfrac12\langle y,Q(x)y\rangle+R(x,y).
\]
Then there is a simple canonical transformation $W\colon D(r/2,s/2)\to D(r,s)$ such that
\begin{equation}\label{eq:Kform}
  H\circ W(\xi,\eta)=a_++\langle\omega,\eta\rangle+\tfrac12\langle\eta,Q_+(\xi)\eta\rangle+R_\ast(\xi,\eta)
\end{equation}
on $D(r/2,s/2)$, with $a_+\in\R$, a symmetric $Q_+\in P_{n\times n}(r/2)$ and $R_\ast\in P(r/2,s/2)$, and
\begin{align}
  \abs{W_\zeta-E}_{r/2,s/2}&\le\Cs\vartheta,\label{eq:Wder}\\
  \abs{W(\xi,0)-(\xi,0)}&\le\Cs\vartheta\,r\qquad\text{for }\xi\in S(r/2),\label{eq:Wzero}\\
  \abs{Q_+-Q}_{r/2}&\le\Cs\vartheta,\label{eq:Qplus}\\
  \abs{R_\ast(\xi,\eta)}&\le\Cs\vartheta\,\frac{\abs\eta^3}{s}\qquad\text{for }(\xi,\eta)\in D(r/2,s/2).\label{eq:Rstar}
\end{align}
\end{theorem}

By \eqref{eq:Kform}, $H\circ W$ is a normal form up to the term $R_\ast$, which vanishes to third order at $\eta=0$. Hence $\{\eta=0\}$ is invariant for $H\circ W$, with linear flow of frequency $\omega$, and its image under $W$ is an invariant torus of $H$.

\begin{remark}
Theorem~\ref{thm:analytic} is essentially \cite[Satz~3.6]{AlbrechtDiss} for $s=r^{\tau+1}$, together with the new estimate \eqref{eq:Wzero}. Compared with the classical formulation (compare \cite{Salamon04}), two features matter for us. First, the constants depend on the perturbation only through $\vartheta=\abs R/s^2$; this makes the powers of $\delta_k$ cancel in Section~\ref{sec:proof}. Second, the case $\tau=n-1$ is admitted; this rests on R\"ussmann's estimate, Lemma~\ref{lem:russmann}. The analytic theorem in \cite{Salamon04} assumes $\tau>n-1$.
\end{remark}

The proof is a Newton iteration. We first solve the linearized equation (Section~\ref{sec:homol}), then construct the transformation for one Newton step as the time-one map of a Hamiltonian flow (Section~\ref{sec:step}), and finally iterate (Section~\ref{sec:iter}). In Sections~\ref{sec:homol} and \ref{sec:step}, $c$ denotes positive constants depending only on $n$, $\tau$, $\gamma$, $\abs C$ and $\abs{C^{-1}}$; their values may change from line to line.

\subsection{The linearized equation}\label{sec:homol}
Let $H=N+R$ with a normal form $N$ and a small $R$. We look for a canonical transformation $Z$, close to the identity, such that $H\circ Z$ is a normal form up to a remainder of order $\abs R^2$. We take $Z$ as the time-one map of the Hamiltonian flow of a function $\Delta S$. To first order, $H\circ Z\approx N+R+\pb N{\Delta S}$, where $\pb FG=\langle F_x,G_y\rangle-\langle F_y,G_x\rangle$ is the Poisson bracket. This leads to the linearized equation
\begin{equation}\label{eq:lin}
  R+\pb N{\Delta S}-\Delta N=0,
\end{equation}
in which $\Delta S$ and $\Delta N$ are unknown, and $\Delta N$ shall have the property that $N+\Delta N$ is again a normal form.

\begin{lemma}\label{lem:homol}
Let $0<4\delta<r\le1$ and $0<s\le\delta^{\tau+1}$. Let $N\in P(r,s)$ be a normal form with $\abs{N_{yy}-C}_{r,s}\le1/(2\abs{C^{-1}})$, and let $R\in P(r,s)$ with $\abs R_{r,s}\le M$. Then there are $\lambda\in\R^n$, $U\in P(r-\delta)$ and $V\in P_n(r-3\delta)$ such that
\[
  \Delta S(x,y):=\langle\lambda,x\rangle+U(x)+\langle V(x),y\rangle
\]
and $\Delta N:=R+\pb N{\Delta S}$ have the following properties: $\Delta N\in P(r-4\delta,s)$, $N+\Delta N$ is a normal form on $D(r-4\delta,s)$, and
\begin{align}
  \abs{\Delta S_x}_{r-4\delta,s}&\le c\,\frac Ms, &
  \abs{\Delta S_y}_{r-3\delta}&=\abs V_{r-3\delta}\le c\,\frac M{s\delta^\tau},\label{eq:Sest}\\
  \abs{\Delta N-\Delta N(x,0)}_{r-4\delta,s/2}&\le c\,M, &
  \abs{\Delta N_{yy}}_{r-4\delta,s/4}&\le c\,\frac M{s^2}.\label{eq:Nest}
\end{align}
\end{lemma}

Here $\Delta N(x,0)$ does not depend on $x$, since $N+\Delta N$ is a normal form. Note that $\Delta S$ is affine in $y$, and that $\Delta S_x=\lambda+U_x+V_x^Ty$ and $\Delta S_y=V$ are periodic, although $\Delta S$ is not.

\begin{proof}
\emph{The conditions on $\Delta S$.} Since $N(x,0)=a$ and $N_y(x,0)=\omega$, we have $N_x(x,0)=0$ and $N_{xy}(x,0)=0$. Inserting $\Delta S_x=\lambda+U_x+V_x^Ty$ and $\Delta S_y=V$, we obtain at $y=0$
\[
  \Delta N(x,0)=R(x,0)-\langle\omega,\lambda\rangle-\partial_\omega U(x).
\]
Differentiating $\Delta N=R+\langle N_x,V\rangle-\langle N_y,\lambda+U_x+V_x^Ty\rangle$ with respect to $y$ and putting $y=0$ gives
\[
  \Delta N_y(x,0)=R_y(x,0)-N_{yy}(x,0)\bigl(\lambda+U_x(x)\bigr)-\partial_\omega V(x).
\]
Hence $N+\Delta N$ is a normal form if and only if
\begin{align}
  \partial_\omega U&=R(\cdot,0)-\mean{R(\cdot,0)},\label{eq:homU}\\
  \partial_\omega V&=R_y(\cdot,0)-N_{yy}(\cdot,0)\,(\lambda+U_x),\label{eq:homV}
\end{align}
where in the first equation we have used that $\partial_\omega U$ has mean value zero; then $\Delta N(x,0)=\mean{R(\cdot,0)}-\langle\omega,\lambda\rangle$.

\emph{Solution of \eqref{eq:homU}.} The right hand side is bounded by $2M$ on $S(r)$ and has mean value zero. Lemma~\ref{lem:russmann} gives $U$ with $\abs U_{r-\delta}\le cM/\delta^\tau$, and Cauchy's estimate gives $\abs{U_x}_{r-2\delta}\le cM/\delta^{\tau+1}\le cM/s$, since $s\le\delta^{\tau+1}$.

\emph{Choice of $\lambda$.} The right hand side of \eqref{eq:homV} has mean value zero if and only if $\mean{N_{yy}(\cdot,0)}\lambda=\mean{R_y(\cdot,0)}-\mean{N_{yy}(\cdot,0)U_x}$. The real matrix $\mean{N_{yy}(\cdot,0)}$ satisfies $\abs{\mean{N_{yy}(\cdot,0)}-C}\le1/(2\abs{C^{-1}})$, so it is invertible with inverse bounded by $2\abs{C^{-1}}$. We define $\lambda\in\R^n$ by this equation. Since $\abs{R_y(x,0)}\le M/s$ and $\abs{N_{yy}}_{r,s}\le\abs C+1/(2\abs{C^{-1}})\le\frac32\abs C$, we get $\abs\lambda\le cM/s$, and the right hand side of \eqref{eq:homV} is bounded by $cM/s$ on $S(r-2\delta)$.

\emph{Solution of \eqref{eq:homV}.} Lemma~\ref{lem:russmann}, applied componentwise on $S(r-2\delta)$, gives $V$ with $\abs V_{r-3\delta}\le cM/(s\delta^\tau)$. Cauchy's estimate gives $\abs{V_x}_{r-4\delta}\le cM/(s\delta^{\tau+1})$, hence $\abs{V_x^Ty}\le n\abs{V_x}\,s\le cM/\delta^{\tau+1}\le cM/s$ on $D(r-4\delta,s)$. Together with the bounds for $\lambda$ and $U_x$, this proves \eqref{eq:Sest}.

\emph{The estimate for $\Delta N$.} Clearly $\Delta N\in P(r-4\delta,s)$. Write $N=a+\langle\omega,y\rangle+\tilde h$; then $\tilde h(x,0)=0$ and $\tilde h_y(x,0)=0$, so Taylor's formula and $\abs{N_{yy}}\le\frac32\abs C$ give $\abs{\tilde h}_{r,s}\le cs^2$. Cauchy's estimate gives $\abs{\tilde h_x}_{r-\delta,s}\le cs^2/\delta$ and $\abs{\tilde h_y}_{r,s/2}\le cs$. Now
\[
  \pb N{\Delta S}=\pb{\tilde h}{\Delta S}-\langle\omega,\Delta S_x\rangle,\qquad
  \langle\omega,\Delta S_x\rangle=\langle\omega,\lambda\rangle+\partial_\omega U+\langle\partial_\omega V,y\rangle,
\]
and with \eqref{eq:homU} this yields
\[
  \Delta N-\Delta N(x,0)=R-R(x,0)+\pb{\tilde h}{\Delta S}-\langle\partial_\omega V,y\rangle .
\]
On $D(r-4\delta,s/2)$, the first difference is bounded by $2M$, and the last term by $n\cdot cM/s\cdot s$, by \eqref{eq:homV}. For the bracket,
\[
  \abs{\pb{\tilde h}{\Delta S}}\le n\abs{\tilde h_x}\abs V+n\abs{\tilde h_y}\abs{\Delta S_x}\le c\,\frac{s^2}{\delta}\,\frac M{s\delta^\tau}+c\,s\,\frac Ms\le cM,
\]
again because $s\le\delta^{\tau+1}$. This proves the first estimate in \eqref{eq:Nest}. Applying Cauchy's estimate twice in $y$, first with radius $s/8$ and then again with radius $s/8$, gives the second.
\end{proof}

\subsection{One Newton step}\label{sec:step}
The flow of a function $\Delta S$ as in Lemma~\ref{lem:homol} is particularly simple: the equations $\dot x=\Delta S_y=V(x)$, $\dot y=-\Delta S_x$ show that $x(t)$ does not depend on the initial value of $y$, and that $y(t)$ depends affinely on it. Hence the time-$t$ maps are simple canonical transformations.

\begin{lemma}\label{lem:flow}
Let $\rho>\varepsilon>0$, $\sigma>0$, $\lambda\in\R^n$, $U\in P(\rho)$, $V\in P_n(\rho)$, and suppose that $G(x,y):=\lambda+U_x(x)+V_x(x)^Ty$ satisfies
\[
  \abs V_\rho\le A\le\frac\varepsilon4,\qquad\abs G_{\rho,\sigma}\le B\le\frac\sigma8 .
\]
Then for $0\le t\le1$ the time-$t$ map $Z_t$ of $\dot x=V(x)$, $\dot y=-G(x,y)$ is defined on $D(\rho-\varepsilon,\sigma/2)$, it is a simple canonical transformation with values in $D(\rho-\varepsilon/2,3\sigma/4)$, and
\[
  \abs{Z_t-\id}_{\rho-\varepsilon,\sigma/2}\le\max\{A,B\},\qquad
  \abs{(Z_t)_\zeta-E}_{\rho-\varepsilon,\sigma/2}\le e^{L}-1,\qquad L:=\frac{2(B+nA)}{\varepsilon}.
\]
\end{lemma}

\begin{proof}
For $\zeta=(\xi,\eta)\in D(\rho-\varepsilon,\sigma/2)$, a solution $(x(t),y(t))$ with initial value $\zeta$ satisfies $\abs{x(t)-\xi}\le At$ and $\abs{y(t)-\eta}\le Bt$ as long as it stays in $D(\rho,\sigma)$. For $t\le1$ it therefore stays in the closed set $\{\abs{\im x}\le\rho-3\varepsilon/4,\ \abs y\le5\sigma/8\}$, which lies in $D(\rho-\varepsilon/2,3\sigma/4)$ at positive distance from its boundary. Since the vector field is periodic in $\re x$, it is bounded and Lipschitz continuous on this set, so the solution can be continued up to $t=1$. The vector field is $2\pi$-periodic in $x$ and real for real arguments; by uniqueness, $Z_t-\id$ is periodic and real. Since the function $\langle\lambda,x\rangle+U(x)+\langle V(x),y\rangle$ generates this flow, $Z_t$ is canonical, and its $x$-component depends only on $\xi$.

The derivative $Z_\zeta(t):=(Z_t)_\zeta(\zeta)$ satisfies $\frac d{dt}Z_\zeta=\left(\begin{smallmatrix}V_x&0\\-G_x&-V_x^T\end{smallmatrix}\right)Z_\zeta$, $Z_\zeta(0)=E$, where the coefficients are evaluated along the solution, that is, in $D(\rho-\varepsilon/2,3\sigma/4)$. There Cauchy's estimate in $x$ gives $\abs{V_x}\le2A/\varepsilon$, $\abs{V_x^T}\le2nA/\varepsilon$ and $\abs{G_x}\le2B/\varepsilon$, so the coefficient matrix is bounded by $L$. Gronwall's inequality gives $\abs{Z_\zeta(t)}\le e^{Lt}$ and then $\abs{Z_\zeta(t)-E}\le\int_0^tLe^{Lt'}\,dt'\le e^L-1$.
\end{proof}

\begin{lemma}[Newton step]\label{lem:step}
There is a constant $c_N\ge1$, depending only on $n$, $\tau$, $\gamma$, $\abs C$ and $\abs{C^{-1}}$, with the following property. Let $0<6\delta<r\le1$ and $0<s\le\delta^{\tau+1}$. Let $H=N+R$, where $N\in P(r,s)$ is a normal form with $\abs{N_{yy}-C}_{r,s}\le1/(2\abs{C^{-1}})$ and $R\in P(r,s)$ satisfies $\abs R_{r,s}\le M\le s^2/c_N$. Then there are a simple canonical transformation $Z\colon D(r-6\delta,s/8)\to D(r-5\delta,s/4)$ and a normal form $N_+\in P(r-6\delta,s/8)$ such that $R_+:=H\circ Z-N_+$ satisfies
\begin{align}
  \abs{R_+}_{r-6\delta,s/8}&\le c_N\frac{M^2}{s^2}, &
  \abs{(N_+)_{yy}-N_{yy}}_{r-6\delta,s/8}&\le c_N\frac{M}{s^2},\label{eq:step1}\\
  \abs{Z_\zeta-E}_{r-6\delta,s/8}&\le c_N\frac{M}{s^2}, &
  \abs{Z-\id}_{r-6\delta,s/8}&\le c_N\frac{M}{s\delta^\tau}.\label{eq:step2}
\end{align}
\end{lemma}

\begin{proof}
Let $\Delta S$ and $\Delta N$ be given by Lemma~\ref{lem:homol}, and denote by $c_1$ the constant there. We apply Lemma~\ref{lem:flow} with $\rho=r-4\delta$, $\varepsilon=2\delta$, $\sigma=s/4$, $A=c_1M/(s\delta^\tau)$ and $B=c_1M/s$. Its hypotheses $A\le\delta/2$ and $B\le s/32$ hold if $c_N\ge32c_1$: then $M\le s^2/c_N$ and $s\le\delta^{\tau+1}$ give $A\le c_1s/(c_N\delta^\tau)\le\delta/32$ and $B\le c_1s/c_N\le s/32$. Let $Z$ be the time-one map, restricted to $D(r-6\delta,s/8)$, and $Z_t$ the time-$t$ maps; they take values in $D(r-5\delta,3s/16)$. Since $\delta\le1$, we have $B\le c_1M/(s\delta^\tau)$, which gives the second estimate in \eqref{eq:step2}. Moreover $L=(B+nA)/\delta\le(n+1)c_1M/(s\delta^{\tau+1})\le(n+1)c_1M/s^2\le\frac12$, and $e^L-1\le2L$; this gives the first estimate in \eqref{eq:step2}.

We put $N_+:=N+\Delta N$, which is a normal form by Lemma~\ref{lem:homol}, and the estimate for $(N_+)_{yy}-N_{yy}=\Delta N_{yy}$ follows from \eqref{eq:Nest}. For a function $F$ which is analytic near the trajectories, $\frac d{dt}F(Z_t(\zeta))=\pb F{\Delta S}(Z_t(\zeta))$. Taylor's formula along the flow therefore gives
\[
  R\circ Z=R+\int_0^1\pb R{\Delta S}\circ Z_t\,dt,\qquad
  N\circ Z=N+\pb N{\Delta S}+\int_0^1(1-t)\,\pb{\pb N{\Delta S}}{\Delta S}\circ Z_t\,dt .
\]
Using $\pb N{\Delta S}=\Delta N-R$ from \eqref{eq:lin}, we obtain
\[
  R_+=R\circ Z+N\circ Z-N-\Delta N=\int_0^1\pb{F_t}{\Delta S}\circ Z_t\,dt,\qquad F_t:=tR+(1-t)\bigl(\Delta N-\Delta N(x,0)\bigr),
\]
where we have used that the constant $\Delta N(x,0)$ does not contribute to the bracket. By \eqref{eq:Nest}, $\abs{F_t}\le(1+c_1)M$ on $D(r-4\delta,s/2)$. By Cauchy's estimate, $\abs{(F_t)_x}\le(1+c_1)M/\delta$ on $D(r-5\delta,s/2)$ and $\abs{(F_t)_y}\le4(1+c_1)M/s$ on $D(r-4\delta,s/4)$. With \eqref{eq:Sest} we get on $D(r-5\delta,s/4)$, which contains the trajectories,
\[
  \abs{\pb{F_t}{\Delta S}}\le n\abs{(F_t)_x}\abs V+n\abs{(F_t)_y}\abs{\Delta S_x}
  \le n(1+c_1)c_1\Bigl(\frac{M^2}{s\delta^{\tau+1}}+\frac{4M^2}{s^2}\Bigr)\le5n(1+c_1)c_1\frac{M^2}{s^2}.
\]
This proves the first estimate in \eqref{eq:step1}, and the lemma holds with $c_N:=\max\{1,\,32c_1,\,2(n+1)c_1,\,5n(1+c_1)c_1\}$. Finally, $R_+\in P(r-6\delta,s/8)$, because $Z-\id$ is periodic and real.
\end{proof}

\subsection{The iteration}\label{sec:iter}
We now prove Theorem~\ref{thm:analytic}. Let $c_N$ be the constant of Lemma~\ref{lem:step}, and define
\begin{equation}\label{eq:cstar}
  t_{\max}:=\min\Bigl\{4^{-4(\tau+1)},\ \frac1{8\abs{C^{-1}}},\ \frac1{20n}\Bigr\},\qquad
  \cs:=\frac{t_{\max}}{c_N\,32^{2(\tau+1)}},\qquad \Cs:=15n\,c_N\,32^{2(\tau+1)} .
\end{equation}
Let $r$, $s=r^{\tau+1}$, $H$, $Q$, $R$ and $\vartheta\le\cs$ be as in the theorem, and put $t_0:=c_N32^{2(\tau+1)}\vartheta\le t_{\max}$.

\emph{The sequences.} For $j\ge0$ let
\[
  \delta_j:=\frac{r}{32}\,4^{-j},\qquad s_j:=\delta_j^{\tau+1},\qquad r_j:=\frac34r+8\delta_j,\qquad t_j:=t_0^{(3/2)^j},\qquad M_j:=\frac{s_j^2\,t_j}{c_N}.
\]
Then $r_0=r$, $r_j\downarrow\frac34r$, $r_j-r_{j+1}=6\delta_j$, $6\delta_j<r_j\le1$, $s_{j+1}=4^{-(\tau+1)}s_j\le s_j/8$, and $s_0=s/32^{\tau+1}\le s$. Since $t_0\le\frac14$, we have $t_{j+1}=t_j\,t_j^{1/2}\le\frac12t_j$, hence $\sum_jt_j\le2t_0$. Moreover $M_0=\vartheta s^2$.

\emph{Induction.} We claim that there are simple canonical transformations $Z_j\colon D(r_{j+1},s_{j+1})\to D(r_j-5\delta_j,s_j/4)$ such that, with $W_0:=\id$ and $W_{j+1}:=W_j\circ Z_j$, we have $H\circ W_j=N_j+R_j$ on $D(r_j,s_j)$, where $N_j$ is a normal form with
\begin{equation}\label{eq:indhyp}
  \abs{R_j}_{r_j,s_j}\le M_j,\qquad\abs{(N_j)_{yy}-Q}_{r_j,s_j}\le\sum_{l<j}t_l,\qquad\abs{(W_j)_\zeta}_{r_j,s_j}\le\prod_{l<j}(1+t_l).
\end{equation}
For $j=0$ put $N_0:=a+\langle\omega,y\rangle+\frac12\langle y,Qy\rangle$ and $R_0:=R$; then \eqref{eq:indhyp} holds because $D(r_0,s_0)\subseteq D(r,s)$. Suppose \eqref{eq:indhyp} holds for $j$. Then $\abs{(N_j)_{yy}-C}\le1/(4\abs{C^{-1}})+2t_0\le1/(2\abs{C^{-1}})$ by \eqref{eq:Qhyp}, and $M_j\le s_j^2/c_N$ since $t_j\le1$. So Lemma~\ref{lem:step} applies with $r=r_j$, $\delta=\delta_j$, $s=s_j$, $M=M_j$. It gives $Z_j$ and $N_{j+1}$, and $R_{j+1}:=H\circ W_{j+1}-N_{j+1}$ satisfies
\[
  \abs{R_{j+1}}_{r_{j+1},s_{j+1}}\le c_N\frac{M_j^2}{s_j^2}=\frac{s_j^2t_j^2}{c_N}\le\frac{s_{j+1}^2\,t_j^{3/2}}{c_N}=M_{j+1};
\]
here we used $t_j^{1/2}\le t_0^{1/2}\le4^{-2(\tau+1)}=(s_{j+1}/s_j)^2$. The remaining estimates follow from \eqref{eq:step1}, \eqref{eq:step2} with $c_NM_j/s_j^2=t_j$, and the chain rule. This completes the induction.

\emph{Estimates for $W_j$.} By \eqref{eq:indhyp}, $\abs{(W_j)_\zeta}_{r_j,s_j}\le e^{2t_0}\le2$. From $(W_{j+1})_\zeta-E=((W_j)_\zeta\circ Z_j-E)(Z_j)_\zeta+((Z_j)_\zeta-E)$ and $\abs{(Z_j)_\zeta-E}\le t_j$ we get by induction $\abs{(W_j)_\zeta-E}_{r_j,s_j}\le\prod_{l<j}(1+t_l)-1\le e^{2t_0}-1\le3t_0$, since $t_0\le4^{-8}$. The segment from $\zeta\in D(r_{j+1},s_{j+1})$ to $Z_j(\zeta)$ lies in the convex set $D(r_j-5\delta_j,s_j/4)\subseteq D(r_j,s_j)$, so the mean value inequality and \eqref{eq:step2} give
\begin{equation}\label{eq:Wdiff}
  \abs{W_{j+1}-W_j}_{r_{j+1},s_{j+1}}\le2\,c_N\frac{M_j}{s_j\delta_j^\tau}=2\,t_j\,\delta_j ,
\end{equation}
since $s_j=\delta_j^{\tau+1}$. This is the place where we keep the factor $\delta_j$, which was not used in \cite{AlbrechtDiss}.

\emph{Convergence on the torus.} All the sets $D(r_j,s_j)$ contain $S(3r/4)\times\{0\}$. By \eqref{eq:Wdiff}, the maps $W_j(\cdot,0)$ converge uniformly on $S(3r/4)$, and
\[
  \abs{W_j(\xi,0)-(\xi,0)}\le\sum_{l\ge0}2t_l\delta_l\le4t_0\delta_0=\frac{t_0}8\,r\qquad\text{for }\xi\in S(3r/4).
\]
Write $W_j=(u_j,w_j)$ and $v_j:=w_j(\cdot,0)$. Then $u_j\to u$ and $v_j\to v$ uniformly on $S(3r/4)$, the limits are analytic, $u-\id$ and $v$ are periodic and real, and $u_{j,\xi}\to u_\xi$, $v_{j,\xi}\to v_\xi$ locally uniformly by Weierstrass' theorem. Since $\abs{(W_j)_\zeta-E}\le3t_0$ at $\eta=0$, we have $\abs{u_\xi-E}_{3r/4}\le3t_0$ and $\abs{v_\xi}_{3r/4}\le3t_0$. Moreover, every $W_j(\xi,0)$ lies in the image of $Z_0$, that is, in $D(r-5\delta_0,s_0/4)$; hence $\abs v_{3r/4}\le s_0/4\le s/4$.

\emph{The transformation $W$.} By Lemma~\ref{lem:simple}, $W_j(\xi,\eta)=(u_j(\xi),v_j(\xi)+u_{j,\xi}(\xi)^{-T}\eta)$, and this formula defines analytic canonical maps on $S(3r/4)\times\C^n$. They converge locally uniformly, together with their derivatives, to
\[
  W(\xi,\eta):=\bigl(u(\xi),\,v(\xi)+u_\xi(\xi)^{-T}\eta\bigr),
\]
which is therefore canonical; we restrict $W$ to $D(r/2,s/2)$. From $\abs{u_\xi-E}\le3t_0\le\frac3{20}$ we get $\abs{u_\xi^{-1}-E}\le4t_0$ and $\abs{u_\xi^{-T}-E}\le4nt_0\le\frac15$ on $S(3r/4)$. For $(\xi,\eta)\in D(r/2,5s/8)$ this gives $\abs{w(\xi,\eta)}\le\abs v+\frac65\abs\eta<\frac{s}{32}+\frac34s<s$, since $\abs v\le s_0/4\le s/32$, and, since $u$ is real for real arguments, $\abs{\im u(\xi)}\le(1+3t_0)\abs{\im\xi}<r$ by the mean value inequality. Hence $W$ maps $D(r/2,5s/8)$ into $D(r,s)$, and in particular $D(r/2,s/2)$ into $D(r,s)$. Clearly $W-\id\in P_{2n}(r/2,s/2)$.

\emph{Estimates \eqref{eq:Wder} and \eqref{eq:Wzero}.} The latter was shown above, since $t_0/8\le\Cs\vartheta$. For the former we write
\[
  W_\zeta-E=\begin{pmatrix}u_\xi-E&0\\v_\xi+\partial_\xi\bigl((u_\xi^{-T}-E)\eta\bigr)&u_\xi^{-T}-E\end{pmatrix}.
\]
For fixed $\eta$ with $\abs\eta<s/2$, the function $\xi\mapsto(u_\xi^{-T}-E)\eta$ is bounded by $2nt_0s$ on $S(3r/4)$, so Cauchy's estimate gives $\abs{\partial_\xi((u_\xi^{-T}-E)\eta)}\le8nt_0s/r\le8nt_0$ on $S(r/2)$. Hence $\abs{W_\zeta-E}_{r/2,s/2}\le3t_0+8nt_0+4nt_0\le15nt_0=\Cs\vartheta$.

\emph{The transformed Hamiltonian.} Let $\mathcal K:=H\circ W$ on $D(r/2,5s/8)$. Since $W_j(\xi,0)\to W(\xi,0)$, $\abs{R_j}_{r_j,s_j}\le M_j\to0$, and $N_j(\xi,0)=:a_j$ is constant, we get
\[
  \mathcal K(\xi,0)=\lim_{j\to\infty}H(W_j(\xi,0))=\lim_{j\to\infty}\bigl(a_j+R_j(\xi,0)\bigr)=\lim_{j\to\infty}a_j=:a_+ .
\]
Next, $(H\circ W_j)_\eta(\xi,0)=(W_j)_\eta^T\,\nabla H(W_j(\xi,0))$ converges to $\mathcal K_\eta(\xi,0)$, and it equals $\omega+(R_j)_\eta(\xi,0)$ with $\abs{(R_j)_\eta(\xi,0)}\le M_j/s_j\to0$. Hence $\mathcal K_\eta(\xi,0)=\omega$. Similarly, since $W_j$ is affine in $\eta$,
\[
  (H\circ W_j)_{\eta\eta}(\xi,0)=(W_j)_\eta^T\,H_{zz}(W_j(\xi,0))\,(W_j)_\eta\ \longrightarrow\ \mathcal K_{\eta\eta}(\xi,0)=:Q_+(\xi).
\] On the other hand it equals $(N_j)_{yy}(\xi,0)+(R_j)_{\eta\eta}(\xi,0)$, where the second term is bounded by $4M_j/s_j^2\to0$, and by \eqref{eq:indhyp} the first term differs from $Q(\xi)$ by at most $\sum_lt_l\le2t_0$. This gives \eqref{eq:Qplus}. Clearly $a_+\in\R$, and $Q_+\in P_{n\times n}(r/2)$ is symmetric.

\emph{The remainder.} Define $R_\ast$ by \eqref{eq:Kform}. By what we have just shown, $R_\ast(\xi,\cdot)=T_3\mathcal K(\xi,\cdot)$ for every $\xi\in S(r/2)$. Now $\mathcal K=N_0\circ W+R\circ W$. The first term is a polynomial of degree at most two in $\eta$, because $N_0$ is quadratic in $y$ and $W(\xi,\eta)$ is affine in $\eta$ with first component independent of $\eta$; so its third order Taylor remainder vanishes. The second term is bounded by $\vartheta s^2$ on $D(r/2,5s/8)$. Lemma~\ref{lem:taylor} with $\sigma=5s/8$ and $\varepsilon=4/5$ gives
\[
  \abs{R_\ast(\xi,\eta)}\le5\vartheta s^2\,\frac{\abs\eta^3}{(5s/8)^3}\le21\,\vartheta\,\frac{\abs\eta^3}{s}\qquad\text{for }\abs\eta\le\frac s2,
\]
which proves \eqref{eq:Rstar}, as $21\le\Cs$. This completes the proof of Theorem~\ref{thm:analytic}.\qed

\section{Proof of Theorem~\ref{thm:main}}\label{sec:proof}

\subsection{The common modulus of continuity}
Put $K(\delta):=K_1(\delta)+\dots+K_{2n}(\delta)+\delta$. Then $K$ is increasing, positive for $\delta>0$, $K(\delta)\to0$ for $\delta\to0$, it majorizes each $K_j$, and $K(N\delta)\le NK(\delta)$ for $N\in\mathbb N$, by the triangle inequality. By \eqref{eq:dini}, $\int_0^1K(\delta)\,\delta^{-1}\,d\delta<\infty$. The summand $\delta$ only serves to make $K$ positive.

\subsection{The scales}
Let $\cs$ and $\Cs$ be the constants of Theorem~\ref{thm:analytic} with $\tau=n-1$; they depend only on $n$, $\gamma$ and $C$. We fix
\[
  q:=\min\bigl\{\tfrac12,\ (2\Cs)^{-1/n}\bigr\},\qquad\text{so that}\qquad\Cs q^n\le\tfrac12,
\]
and, with $\delta_0\in(0,q]$ to be chosen below,
\[
  \delta_k:=q^k\delta_0,\qquad s_k:=\delta_k^{\,n}=\delta_k^{\tau+1},\qquad\varrho_k:=\delta_k/q.
\]
Then $\varrho_0\le1$, and $D(\delta_{k+1},s_{k+1})\subseteq D(\delta_k/2,s_k/2)$ since $q\le\frac12$ and $q^n\le\frac12$. Theorem~\ref{thm:approx}, applied with $\nu=2n$, $m=2n$, $f=h$ and the sequence $(\varrho_k)$, gives analytic functions $h_k$ on $S'(\varrho_k)$. Let
\[
  H_k(x,y):=a+\langle\omega,y\rangle+\tfrac12\langle y,Cy\rangle+h_k(x,y).
\]
Since $D(\delta,s)\subseteq S'(\delta)$ for $s\le\delta$, the function $H_k$ is defined on $D(\varrho_k,s)$ for $s\le\varrho_k$. Using $K(\delta_k/q)\le(1/q+1)K(\delta_k)\le(2/q)K(\delta_k)$ and $\delta_k^{2n}=q^{-2n}s_{k+1}^2$, we obtain
\begin{equation}\label{eq:hdiff}
  \abs{h_{k+1}-h_k}_{S'(\delta_k)}\le\ca K(\varrho_k)\varrho_k^{2n}\le b_k\,s_{k+1}^2,\qquad b_k:=2\ca\,q^{-4n-1}K(\delta_k),
\end{equation}
where $S'(\delta_k)=S'(\varrho_{k+1})$.

\subsection{Where the Dini condition enters}
We shall apply Theorem~\ref{thm:analytic} in step $k$ with $r=\delta_k$ and $s=s_k$, and with $\vartheta=\vartheta_k$ defined by
\begin{equation}\label{eq:theta}
  \vartheta_0:=b_0,\qquad\vartheta_{k+1}:=\tfrac12\vartheta_k+b_k\qquad(k\ge0).
\end{equation}
The first summand accounts for the cubic remainder $R_\ast$ of the previous step, which shrinks with the factor $\Cs q^n\le\frac12$ when we pass from $s_k$ to $s_{k+1}$; the second accounts for the new correction $h_{k+1}-h_k$. Summing \eqref{eq:theta} gives $\sum_k\vartheta_k\le2\vartheta_0+2\sum_kb_k\le4\sum_kb_k$.

\begin{lemma}\label{lem:dini}
$\delta_0\in(0,q]$ can be chosen, depending only on $n$, $\gamma$, $C$ and $K$, such that
\begin{equation}\label{eq:thetasum}
  \sum_{k=0}^\infty\vartheta_k\le\Theta:=\min\Bigl\{\cs,\ \frac{1}{4\Cs\abs{C^{-1}}},\ \frac{\log(3/2)}{\Cs}\Bigr\}.
\end{equation}
\end{lemma}

\begin{proof}
Since $K$ is increasing, $K(\delta_k)\log(1/q)\le\int_{\delta_k}^{\delta_{k-1}}K(t)\,t^{-1}\,dt$ for $k\ge1$, so
\[
  \sum_{k=0}^\infty\vartheta_k\le\frac{8\ca}{q^{4n+1}}\Bigl(K(\delta_0)+\frac1{\log(1/q)}\int_0^{\delta_0}\frac{K(t)}t\,dt\Bigr),
\]
and the right hand side tends to zero as $\delta_0\to0$.
\end{proof}

This is the only place where the Dini condition is used. We fix $\delta_0$ as in Lemma~\ref{lem:dini} and put
\[
  \eA:=\frac{\vartheta_0\,s_0^2}{\ca}=2q^{-4n-1}K(\delta_0)\,\delta_0^{2n}>0 .
\]

\subsection{The iteration}
\begin{proposition}\label{prop:iter}
Let $\abs h_{\R^{2n}}\le\eA$. Then for every $k\ge0$ there is a simple canonical transformation $Z_k$ on $D(\delta_k/2,s_k/2)$ with values in $S'(\delta_k)$ such that
\begin{equation}\label{eq:normal}
  H_k\circ Z_k(\xi,\eta)=a_{k+1}+\langle\omega,\eta\rangle+\tfrac12\langle\eta,Q_{k+1}(\xi)\eta\rangle+R^\ast_k(\xi,\eta)
\end{equation}
on $D(\delta_k/2,s_k/2)$, with $a_{k+1}\in\R$, a symmetric $Q_{k+1}\in P_{n\times n}(\delta_k/2)$, and
\begin{align}
  \abs{(Z_k)_\zeta-E}_{\delta_k/2,s_k/2}&\le\exp\Bigl(\Cs\sum_{l\le k}\vartheta_l\Bigr)-1\le\frac12,\label{eq:p1}\\
  \abs{Q_{k+1}-C}_{\delta_k/2}&\le\Cs\sum_{l\le k}\vartheta_l\le\frac{1}{4\abs{C^{-1}}},\label{eq:p2}\\
  \abs{R^\ast_k(\xi,\eta)}&\le\Cs\vartheta_k\,\frac{\abs\eta^3}{s_k}\qquad\text{on }D(\delta_k/2,s_k/2),\label{eq:p3}\\
  \abs{Z_k(\xi,0)-Z_{k-1}(\xi,0)}&\le\tfrac32\Cs\vartheta_k\,\delta_k\qquad\text{for }\xi\in S(\delta_k/2),\label{eq:p4}
\end{align}
where $Z_{-1}:=\id$.
\end{proposition}

\begin{proof}
The second inequalities in \eqref{eq:p1} and \eqref{eq:p2} follow from \eqref{eq:thetasum}. We first note a consequence of \eqref{eq:p1}: if $\zeta\in D(\delta_k/2,s_k/2)$, then $\re\zeta$ lies in the same convex set and $Z_k(\re\zeta)$ is real, so $\abs{\im Z_k(\zeta)}\le\abs{Z_k(\zeta)-Z_k(\re\zeta)}\le\frac32\abs{\im\zeta}<\frac34\delta_k$, because $\abs{\im\eta}\le\abs\eta<s_k/2\le\delta_k/2$. Hence $Z_k$ maps $D(\delta_k/2,s_k/2)$ into $S'(\delta_k)$, once \eqref{eq:p1} is known.

\emph{The case $k=0$.} On $D(\delta_0,s_0)\subseteq S'(\varrho_0)$ we have $H_0=a+\langle\omega,y\rangle+\frac12\langle y,Cy\rangle+h_0$ with $\abs{h_0}\le\ca\abs h_{\R^{2n}}\le\ca\eA=\vartheta_0s_0^2$. Theorem~\ref{thm:analytic} with $r=\delta_0$, $Q=C$, $R=h_0$ and $\vartheta=\vartheta_0\le\cs$ gives $W$, and we put $Z_0:=W$. Then \eqref{eq:normal} and \eqref{eq:p2}--\eqref{eq:p4} follow from \eqref{eq:Kform}--\eqref{eq:Rstar}, and \eqref{eq:p1} from \eqref{eq:Wder} and $1+t\le e^t$.

\emph{From $k$ to $k+1$.} Since $Z_k$ maps $D(\delta_k/2,s_k/2)$ into $S'(\delta_k)=S'(\varrho_{k+1})$, where $h_{k+1}$ is defined, we may write on $D(\delta_{k+1},s_{k+1})$
\[
  H_{k+1}\circ Z_k=a_{k+1}+\langle\omega,\eta\rangle+\tfrac12\langle\eta,Q_{k+1}\eta\rangle+R,\qquad R:=R^\ast_k+(h_{k+1}-h_k)\circ Z_k .
\]
The function $R$ is periodic in $\xi$, because $Z_k-\id$ is, and real for real arguments. By \eqref{eq:p3}, $\abs{R^\ast_k}\le\Cs\vartheta_ks_{k+1}^3/s_k=\Cs q^n\vartheta_ks_{k+1}^2\le\frac12\vartheta_ks_{k+1}^2$ on $D(\delta_{k+1},s_{k+1})$, and by \eqref{eq:hdiff}, $\abs{(h_{k+1}-h_k)\circ Z_k}\le b_ks_{k+1}^2$. Hence $\abs R_{\delta_{k+1},s_{k+1}}\le\vartheta_{k+1}s_{k+1}^2$ with $\vartheta_{k+1}\le\cs$, and $Q_{k+1}$ satisfies \eqref{eq:Qhyp} on $S(\delta_{k+1})$ by \eqref{eq:p2}. Theorem~\ref{thm:analytic} with $r=\delta_{k+1}$, $s=s_{k+1}$ and $\vartheta=\vartheta_{k+1}$ gives $W\colon D(\delta_{k+1}/2,s_{k+1}/2)\to D(\delta_{k+1},s_{k+1})$, and we put $Z_{k+1}:=Z_k\circ W$. Then $H_{k+1}\circ Z_{k+1}=(H_{k+1}\circ Z_k)\circ W$ has the form \eqref{eq:normal} by \eqref{eq:Kform}, and \eqref{eq:p2}, \eqref{eq:p3} for $k+1$ follow from \eqref{eq:Qplus} and \eqref{eq:Rstar}. For \eqref{eq:p1} we write $(Z_{k+1})_\zeta-E=((Z_k)_\zeta\circ W-E)W_\zeta+(W_\zeta-E)$; with $\alpha:=\exp(\Cs\sum_{l\le k}\vartheta_l)-1$ and $\beta:=\Cs\vartheta_{k+1}$ this gives the bound $\alpha(1+\beta)+\beta\le(1+\alpha)e^\beta-1$. Finally, $W(\xi,0)$ and $(\xi,0)$ lie in the convex set $D(\delta_{k+1},s_{k+1})\subseteq D(\delta_k/2,s_k/2)$, on which $\abs{(Z_k)_\zeta}\le\frac32$, so \eqref{eq:Wzero} gives
\[
  \abs{Z_{k+1}(\xi,0)-Z_k(\xi,0)}=\abs{Z_k(W(\xi,0))-Z_k(\xi,0)}\le\tfrac32\Cs\vartheta_{k+1}\delta_{k+1}.\qedhere
\]
\end{proof}

\subsection{Convergence}
\begin{proof}[Proof of Theorem~\ref{thm:main}]
Write $Z_k(\xi,0)=(X_k(\xi),Y_k(\xi))$. By \eqref{eq:p4} and Lemma~\ref{lem:dini}, the maps $Z_k(\cdot,0)$ converge uniformly on $\R^n$ to a continuous map $(X,Y)$. To control the derivatives, let $\xi\in\R^n$. The function $\xi'\mapsto Z_k(\xi',0)-Z_{k-1}(\xi',0)$ is analytic on $S(\delta_k/2)$, which contains the ball of radius $\delta_k/2$ around $\xi$. Cauchy's estimate and \eqref{eq:p4} give
\[
  \bigl|\partial_\xi\bigl(Z_k(\xi,0)-Z_{k-1}(\xi,0)\bigr)\bigr|\le\frac{\frac32\Cs\vartheta_k\delta_k}{\delta_k/2}=3\Cs\vartheta_k .
\]
By Lemma~\ref{lem:dini} the right hand side is summable. Hence the derivatives converge uniformly on $\R^n$ as well, $X$ and $Y$ are continuously differentiable, $X-\id$ and $Y$ are periodic, and by \eqref{eq:p1}, $\abs{X_\xi-E}\le\frac12$ and $\abs{Y_\xi}\le\frac12$ on $\R^n$. Thus $X-\id$ is Lipschitz continuous with constant $\frac12$; so $X$ is injective, it is surjective since $\xi\mapsto x-(X(\xi)-\xi)$ is a contraction for every $x$, and $X^{-1}$ is $C^1$ by the inverse function theorem. Moreover $X^{-1}-\id$ is periodic.

\emph{Solutions.} By \eqref{eq:normal}, $(H_k\circ Z_k)_\xi(\xi,0)=0$ and $(H_k\circ Z_k)_\eta(\xi,0)=\omega$, since the quadratic and the cubic terms do not contribute at $\eta=0$. Hence $\zeta(t)=(\xi_0+\omega t,0)$ solves the Hamiltonian system for $H_k\circ Z_k$, and since $Z_k$ is canonical, $z_k(t):=Z_k(\xi_0+\omega t,0)$ solves $\dot z=J\nabla H_k(z)$, that is, $z_k(t)=z_k(0)+\int_0^tJ\nabla H_k(z_k(\sigma))\,d\sigma$. As $k\to\infty$, $z_k\to z:=(X,Y)(\xi_0+\omega\,\cdot)$ uniformly, and $\nabla H_k\to\nabla H$ uniformly on $\R^{2n}$ by Theorem~\ref{thm:approx}. Since $\nabla H$ is uniformly continuous on bounded sets, $z(t)=z(0)+\int_0^tJ\nabla H(z(\sigma))\,d\sigma$, which proves (1).

\emph{The torus.} Let $g:=Y\circ X^{-1}$; it is $C^1$ and periodic. By Lemma~\ref{lem:simple}, $X_{k,\xi}^TY_{k,\xi}$ is symmetric for every $k$, and in the limit $X_\xi^TY_\xi$ is symmetric. Since $g_x\circ X=Y_\xi X_\xi^{-1}$, also $g_x=X_\xi^{-T}(X_\xi^TY_\xi)X_\xi^{-1}$ is symmetric. Hence $g=\nabla\Phi$ for a function $\Phi\in C^2(\R^n)$. Next, $H_k(Z_k(\xi,0))=a_{k+1}$ by \eqref{eq:normal}; the left hand side converges to $H(X(\xi),Y(\xi))$, so $a_{k+1}$ converges and $H(x,\nabla\Phi(x))$ is constant. The torus $\{(x,\nabla\Phi(x))\}=\{(X(\xi),Y(\xi)):\xi\in\R^n\}$ is the graph of a $C^1$ map, hence a $C^1$ submanifold. By (1) it is invariant, and every orbit on it is dense, since the components of $\omega$ are rationally independent by \eqref{eq:dioph}.
\end{proof}

\begin{remark}
The same proof works for an integer $\tau\ge n-1$, if $h$ is of class $C^{2\tau+2}$ and the pure derivatives of order $2\tau+2$ satisfy the Dini condition. One takes $s_k=\delta_k^{\tau+1}$, chooses $q$ with $\Cs q^{\tau+1}\le\frac12$, and applies Theorem~\ref{thm:approx} with $m=2\tau+2$; the powers of $\delta_k$ again cancel. For $\tau=n$ this gives a result for almost all $\omega$, with two more derivatives.
\end{remark}

\subsection*{Acknowledgements}
This paper was prepared with substantial assistance from the AI system Claude (Anthropic; model Claude Opus 5.5). It was used to restructure and translate the material of \cite{AlbrechtDiss}, to draft the proofs, and to check the estimates. The author has verified all arguments and references and takes full responsibility for the content. Any remaining errors are the author's own.


\begin{thebibliography}{99}\small
\bibitem{Achieser67} N.\,I. Achieser, \emph{Vorlesungen \"uber Approximationstheorie}, 2nd ed., Akademie-Verlag, Berlin, 1967.
\bibitem{AlbrechtDiss} J. Albrecht, \emph{\"Uber die Stabilit\"at invarianter Tori in Hamiltonschen Systemen, die bis auf eine endlich oft differenzierbare St\"orung analytisch und integrabel sind}, Dissertation, Johannes Gutenberg-Universit\"at Mainz, 2005, doi:10.25358/openscience-1904.
\bibitem{Albrecht05F} J. Albrecht, Kleine Nenner bei $2\pi$-periodischen Vektorfeldern, in: H.\,P. Delfosse, H.\,R. Yousefi (eds.), \emph{Festschrift f\"ur Michael Albrecht zum 65. Geburtstag}, Traugott Bautz, Nordhausen, 2005, 371--388.
\bibitem{Albrecht07} J. Albrecht, On the existence of invariant tori in nearly-integrable Hamiltonian systems with finitely differentiable perturbations, \emph{Regul. Chaotic Dyn.} \textbf{12} (2007), 281--320.
\bibitem{Arnold63} V.\,I. Arnold, Proof of a theorem of A.\,N. Kolmogorov on the invariance of quasi-periodic motions under small perturbations of the Hamiltonian, \emph{Russian Math. Surveys} \textbf{18} (1963), no.~5, 9--36.
\bibitem{Arnold89} V.\,I. Arnold, \emph{Mathematical methods of classical mechanics}, 2nd ed., Springer, New York, 1989.
\bibitem{Bounemoura20} A. Bounemoura, Positive measure of KAM tori for finitely differentiable Hamiltonians, \emph{J. \'Ec. polytech. Math.} \textbf{7} (2020), 1113--1132.
\bibitem{ChengWang13} C.-Q. Cheng, L. Wang, Destruction of Lagrangian torus for positive definite Hamiltonian systems, \emph{Geom. Funct. Anal.} \textbf{23} (2013), 848--866.
\bibitem{Herman83} M.\,R. Herman, Sur les courbes invariantes par les diff\'eomorphismes de l'anneau, Vol.~1, \emph{Ast\'erisque} \textbf{103--104} (1983).
\bibitem{Kolmogorov54} A.\,N. Kolmogorov, On conservation of conditionally periodic motions for a small change in Hamilton's function, \emph{Dokl. Akad. Nauk SSSR} \textbf{98} (1954), 527--530.
\bibitem{Koudjinan20} C.\,E. Koudjinan, A KAM theorem for finitely differentiable Hamiltonian systems, \emph{J. Differential Equations} \textbf{269} (2020), 4720--4750.
\bibitem{Moser62} J. Moser, On invariant curves of area-preserving mappings of an annulus, \emph{Nachr. Akad. Wiss. G\"ottingen Math.-Phys. Kl. II} (1962), 1--20.
\bibitem{Moser70} J. Moser, On the construction of almost periodic solutions for ordinary differential equations, in: \emph{Proc. Internat. Conf. on Functional Analysis and Related Topics (Tokyo, 1969)}, Univ. of Tokyo Press, Tokyo, 1970, 60--67.
\bibitem{Moser73} J. Moser, \emph{Stable and random motions in dynamical systems}, Princeton University Press, Princeton, 1973.
\bibitem{Paulusch26} J. Paulusch, Linearization of vector fields on the $n$-torus of class $C^n$ with Dini continuous highest derivatives, preprint, 2026.
\bibitem{Poschel80} J. P\"oschel, \"Uber invariante Tori in differenzierbaren Hamiltonschen Systemen, \emph{Bonner Math. Schriften} \textbf{120} (1980).
\bibitem{Poschel82} J. P\"oschel, Integrability of Hamiltonian systems on Cantor sets, \emph{Comm. Pure Appl. Math.} \textbf{35} (1982), 653--696.
\bibitem{Poschel21} J. P\"oschel, KAM below $C^n$, preprint, 2021, arXiv:2104.01866.
\bibitem{Russmann70} H. R\"ussmann, Kleine Nenner I: \"Uber invariante Kurven differenzierbarer Abbildungen eines Kreisringes, \emph{Nachr. Akad. Wiss. G\"ottingen Math.-Phys. Kl. II} (1970), 67--105.
\bibitem{Russmann75} H. R\"ussmann, On optimal estimates for the solutions of linear partial differential equations of first order with constant coefficients on the torus, in: J. Moser (ed.), \emph{Dynamical systems, theory and applications}, Lecture Notes in Physics \textbf{38}, Springer, Berlin, 1975, 598--624.
\bibitem{Russmann79} H. R\"ussmann, Konvergente Reihenentwicklungen in der St\"orungstheorie der Himmelsmechanik, in: K. Jacobs (ed.), \emph{Selecta Mathematica V}, Springer, Berlin, 1979, 93--260.
\bibitem{Russmann83} H. R\"ussmann, On the existence of invariant curves of twist mappings of an annulus, in: J. Palis (ed.), \emph{Geometric dynamics}, Lecture Notes in Math. \textbf{1007}, Springer, Berlin, 1983, 677--718.
\bibitem{Salamon04} D.\,A. Salamon, The Kolmogorov--Arnold--Moser theorem, \emph{Math. Phys. Electron. J.} \textbf{10} (2004), Paper 3.
\bibitem{TongLi25} Z. Tong, Y. Li, Towards continuity: universal frequency-preserving KAM persistence and remaining regularity in nearly non-integrable Hamiltonian systems, \emph{Commun. Contemp. Math.} \textbf{27} (2025), doi:10.1142/S021919972450038X.
\bibitem{Wang23} L. Wang, Quantitative destruction and persistence of Lagrangian torus in Hamiltonian systems, preprint, 2023, arXiv:2312.01695.
\bibitem{Zehnder75} E. Zehnder, Generalized implicit function theorems with applications to some small divisor problems I, \emph{Comm. Pure Appl. Math.} \textbf{28} (1975), 91--140.
\bibitem{Zehnder76} E. Zehnder, Generalized implicit function theorems with applications to some small divisor problems II, \emph{Comm. Pure Appl. Math.} \textbf{29} (1976), 49--111.
\end{thebibliography}
\end{document}